\documentclass{amsart}[12pt]

\usepackage[cmtip,all]{xy}
\usepackage[utf8]{inputenc}
\usepackage{amsmath}
\usepackage{mathrsfs}
\usepackage{amssymb}
\usepackage{mathtools}
\usepackage{url}
\usepackage[top=1.3in, bottom=1.3in, left=1.3in, right=1.3in]{geometry}
\usepackage{pxfonts}
\usepackage{tikz-cd}
\usetikzlibrary{matrix,arrows}
\usepackage{hyperref}

\newtheorem{theorem}{Theorem}[section]
\newtheorem{lemma}[theorem]{Lemma}
\newtheorem{cor}[theorem]{Corollary}
\newtheorem{prop}[theorem]{Proposition}
\theoremstyle{definition}
\newtheorem{defn}[theorem]{Definition}

\newtheorem{claim}[theorem]{Claim}

\theoremstyle{remark}
\newtheorem{rmk}[theorem]{Remark}

\numberwithin{equation}{section}

\newcommand{\abs}[1]{\lvert#1\rvert}

\newcommand*{\sheafhom}{\mathscr{H}\kern -.5pt om}

\begin{document}

\title[Counterexamples of Kodaira vanishing]{Counterexamples of Kodaira vanishing for smooth surfaces of general type in positive characteristic}%
\author{Xudong Zheng}%
\address{University of Illinois at Chicago, 851 South Morgan Street, Science and Engineering Offices m/c249, Chicago, IL 60607 USA}%
\email{xzheng20@uic.edu}%

\maketitle
\date{\today}%

\begin{abstract}
We generalize the construction of Raynaud \cite{R78} of smooth projective surfaces of general type in positive characteristic that violate the Kodaira vanishing theorem. This corrects an earlier paper \cite{T10} of the same purpose. These examples are smooth surfaces fibered over a smooth curve whose direct images of the relative dualizing sheaves are not nef, and they violate Koll\'ar's vanishing theorem. Further pathologies on these examples include the existence of non-trivial vector fields and that of non-closed global differential 1-forms.
\end{abstract}

\maketitle

\section{Introduction}

The classical Kodaira vanishing theorem on complex projective manifolds asserts the vanishing of cohomology groups of anti-ample line bundles below the dimension. This is well known to be false in positive characteristic. Among the examples, Raynaud constructed for any prime number $p$ a smooth surface of characteristic $p$ with an ample line bundle violating the Kodaira vanishing (see \cite{R78}). Several authors later generalized Raynaud's examples, which all exhibit certain geometric properties only possible in characteristic $p$, such as singular fibration structure and uniruledness for surfaces of general type (see \cite{E88, T92, Mu13}). In fact, the existence of a fibration over a smooth curve whose generic fiber is singular can characterize such Kodaira non-vanishing examples (see \cite{Mu13}). More recently in \cite{T10}, the author claims discovery of similar examples. The idea is to use degree-prime-to-$p$ ramified covers of certain ruled surfaces over a \textit{Tango curve} of sufficiently high genus. However, some of the results are erroneous.

The purpose of this note is to correct a lemma (\cite[Lemma 4]{T10}) of central importance in \cite{T10} and in turn to construct some new examples of Kodaira non-vanishing for smooth surfaces. In light of the few characteristic-free results on pluricanonical systems and rank 2 vector bundles on surfaces not of general type (see \cite{E88, ShB91a}), the current focus is on surfaces of general type. Yet little is known on the classification of general type surfaces with certain positive characteristic pathologies (see \cite{L08, LS09}). For this matter, we would like to further understand the known examples of Kodaira non-vanishing and to come up with new ones. It turns out that our examples also possess some other well-known characteristic $p$ pathologies.

Below we first quickly characterize the known examples in terms of the positivity of the ample line bundle. There are mainly two types of Kodaira non-vanishing examples of surfaces. Suppose $X$ is a smooth projective surface of general type over an algebraically closed field $k$ of characteristic $p$, and $L$ is an ample line bundle on $X$ such that $H^1(X, L^{-1}) \neq 0$.

\textbf{Type 1}. When $L^{p - 1} \otimes \omega_X^{-1}$ is pseudo-effective, i.e., for any nef line bundle $M$ on $X$ we have $(M \cdot L^{p - 1} \otimes \omega_X^{-1}) > 0$, then $X$ is purely inseparably uniruled. One can construct the ruled cover $f: Y \to X$ purely inseparable of degree $p$ where $Y$ is normal projective with precisely $\omega_Y = f^*(L^{1 - p} \otimes \omega_X)$. The fact that there are $K_Y$-negative curves would produce rational curves that are rulings. This type is special in the sense of the following result of Shepherd-Barron:

\begin{prop}\cite[Prop. 21]{ShB91a}\label{regular}
If $L^{p - 1} \otimes K_X^{-1}$ is nef and big, and $H^1(X, L^{-1}) \neq 0$. Then $q(X) = 0$, i.e., $X$ is regular.
\end{prop}

Indeed, Raynaud's examples and the construction of this article fall into the second type.

\textbf{Type 2}. When $L^{p - 1} \otimes \omega_X^{-1}$ is not pseudo-effective. In this case, the well-known Bend-and-Break technique to construct rational curves on varieties no longer works. There still exists some purely inseparable cover $Y$, but $Y$ could be of general type.

We briefly describe the construction and the non-vanishing results. Suppose $(C, L)$ is a Tango curve of genus $g$ (see Definition \ref{tango}) with $2g - 2$ divisible by $p = \mathrm{char}(k)$, where $L \in \mathrm{Pic}(C)$ is an ample line bundle of degree $l \coloneqq \deg(L) \mid \mathrm{gcd}(\frac{2g - 2}{p}, p + 1)$ such that $H^1(C, L^{-1}) \neq 0$, and write $m = (p + 1)/l$. Suppose $\mathcal{E}$ is a rank 2 vector bundle corresponding to a non-zero class $0 \neq \eta \in H^1(C, L^{-1})$ and $P \coloneqq \mathbb{P}(\mathcal{E})$ is the associated ruled surface over $C$. That the curve $C$ being Tango gives rise to a section $E \subset P$ and a \textit{multi-section} $C''$ of $\pi: P \to C$, which are disjoint integral curves on $P$. Let $N$ be a line bundle on $C$ such that $N^l = L$, and let $\psi: X \to P$ be the degree $l$ cover ramified along $E + C''$ given by the line bundle
\[
M \coloneqq \mathcal{O}_P(- m) \otimes \pi^*N^{p}.
\]

We write $\phi: X \xrightarrow{\psi} P \xrightarrow{\pi} C$ for the composition fibration. Finally, let $Z = \mathcal{O}_X(\tilde{E}) \otimes \phi^*N$ be a line bundle on $X$, where $\tilde{E}$ is the reduced preimage of $E$ in $X$. Then $Z$ is ample, and for any positive integer $n \leq \lfloor \dfrac{l}{2} \rfloor$, we prove that $H^1(X, Z^{-n}) \neq 0$ (see Theorem \ref{nonvan1}). Moreover, for any positive integers $a$ and $b$, the line bundle $Z_{a, b} \coloneqq \mathcal{O}_X(a\tilde{E}) \otimes \phi^*N^b$ on $X$ is ample. For any pairs of positive integers $a$ and $b$ such that $a \leq l - 1$ and $b \leq l - a$, we have $H^1(X, Z^{-1}_{a, b}) \neq 0$ (see Theorem \ref{nonvan2}).

Besides Kodaira vanishing, a few other characteristic zero results also fail on $X$. In characteristic zero, Fujita shows that $f_*\omega_{X/C}$ is nef for any 2-dimensional K\"ahler fibered space $f: X \to C$ over a smooth curve (see \cite[Main Thm.]{F78}). Szpiro shows that for a non-isotrivial semi-stable fibration with generic fiber smooth of genus at least 2 the relative dualizing sheaf is nef in arbitrary characteristic (see \cite[Th\'eor\`eme 1]{S79}).  On the other hand, Moret-Bailly constructs an example of a family of genus 2 curves over $\mathbb{P}^1$ where the direct image of the relative dualizing sheaf is not nef (see \cite{MB81}). The current examples indicate that generically singular fibration could be an obstruction to the nefness of $\phi_*\omega_{X/C}$ (see Corollary \ref{nonnef}). 

Moreover, the same computation immediately shows that Koll\'ar's vanishing theorem does not hold true for these surfaces. Briefly, Koll\'ar's vanishing theorem states that for a projective surjective morphism of projective varieties $\pi: X \to Y$ with $X$ smooth and $Y$ normal over $\mathbb{C}$ and for any ample line bundle $A$ on $Y$ the cohomologies $H^i(Y, \mathbf{R}^j\pi_*\omega_X \otimes A)$ vanish for all $i + j > 0$. However, it turns out that $H^1(C, \phi_*\omega_X \otimes N) \neq 0$ in our example (see Corollary \ref{kollar}).

The paper is organized as follows. In section 2, we recall the construction of Tango and Raynaud with all the proofs omitted, and we point out the errors in \cite{T10}. In section 3, we show the corresponding corrections without too much repetition of the proofs existing in the literature. \cite[Lemma 4]{T10} is corrected by Proposition \ref{correctlemma} and \cite[Thm. 18]{T10} is corrected by Lemma \ref{directimage}. Using spectral sequence argument, we provide the examples in Theorem \ref{nonvan1} and \ref{nonvan2} as corrections to \cite[Thm. 19 and Thm. 21]{T10} respectively. As a corollary, we show that our examples are fibrations such that the direct image of the relative dualizing sheaf is not nef (see Corollary \ref{nonnef}), and examples which violate the Koll\'ar's vanishing theorem (see Corollary \ref{kollar}). Finally, the last section collects some ideas behind the construction that are well-known to experts and states some generalizations of results that are at most implicit in the literature. For example, we prove an effective vanishing for ample line bundles on surfaces not of general type (see Corollary \ref{effective}). Finally our examples are also surfaces of general type with non-trivial vector fields, on which not every global differential 1-form is closed (see Section \ref{globalvec}).

\textbf{Acknowledgements}. The author is grateful to Professors Lawrence Ein and Kevin Tucker for their many useful suggestions, and to Professors Christian Liedtke and Adrian Langer for the inspiration to include the last section of the paper. He thanks Professor John Lesieutre for a careful inspection of the previous draft. He also thanks Professors Wenfei Liu, Lance Miller, Mircea Musta\c{t}\u{a}, Mihnea Popa, and Shunsuke Takagi for reading earlier versions of the paper.

\section{Preliminaries}

This section is a brief review of (generalized) Raynaud's examples where the details are left for the reader to consult to the original papers. Throughout this note, we do not distinguish tensor product of line bundles and additive operation of Cartier divisors, and we would use the notation $L^{\otimes 2}$ to eliminate ambiguity from self-intersection of the line bundle $L$. Schemes are over an algebraically closed field $k$ of characteristic $p > 0$. For a scheme $X$ we denote by $F = F_X: X \to X$ the absolute Frobenius endomorphism of $X$.

\subsection{Tango curves}
Let $C$ be a smooth projective curve of genus $g \geq 2$, and write $K(C)$ for its function field. Denote by $K(C)^p = \{f^p \mid f \in K(C)\}$ the subfield of $p$-th powers.

\begin{defn}\label{tango}\cite{T72}
A smooth projective curve $C$ of genus $g \geq 2$ is a \textit{Tango curve} if there is an ample effective divisor $D$ on $C$ such that $(df) = pD$ for some $f \in K(C) \setminus K(C)^p$. The associated line bundle $L = \mathcal{O}_C(D)$ is called a \textit{Tango structure} on $C$.
\end{defn}

Let $C$ be a Tango curve with a Tango structure $L = \mathcal{O}_C(D)$, and consider the sequence:
\begin{equation}\label{4term}
0 \to \mathcal{O}_C \to F_*\mathcal{O}_C \to F_*\Omega^1_C \xrightarrow{c} \Omega^1_C \to 0,
\end{equation}
where the rightmost $\mathcal{O}_C$-module map $c$ is the Cartier operator. Write $B^1 \coloneqq \mathrm{ker}(c)$, then $B^1 = \mathrm{coker}[\mathcal{O}_C \to F_*\mathcal{O}_C]$, i.e., \ref{4term} splits into the following two exact sequences of $\mathcal{O}_C$-modules:
\begin{align}\label{b1}
& 0 \to \mathcal{O}_C \to F_*\mathcal{O}_C \to B^1 \to 0\\
& 0 \to B^1 \to F_*\Omega^1_C \xrightarrow{c} \Omega^1_C \to 0.
\end{align}
Tensoring with $L^{-1}$ and taking cohomology we have

\begin{equation}\label{les1}
0 \to H^0(C, L^{-1}) \to H^0(C, F_*\mathcal{O}_C \otimes L^{-1}) \to H^0(C, B^1 \otimes L^{-1}) \xrightarrow{\beta} H^1(C, L^{-1}) \xrightarrow{\gamma} H^1(C, F_*\mathcal{O}_C \otimes L^{-1})
\end{equation}
and
\begin{equation}\label{les2}
0 \to H^0(C, B^1 \otimes L^{-1}) \to H^0(C, F_*\Omega^1_C \otimes L^{-1}) \xrightarrow{\alpha} H^0(C, \Omega_C^1\otimes L^{-1})
\end{equation}

Here in \ref{les2} $\mathrm{ker}(\alpha) = H^0(C, B^1 \otimes L^{-1}) = \{df \mid f \in K(C), {df} \geq pD\} \neq 0$ since $C$ is Tango. On the other hand, looking at \ref{les1} and suppose that there is a section $0 \neq s \in H^0(C, B^1 \otimes L^{-1})$ such that $0 \neq \beta(s) \in H^1(C, L^{-1})$, then $\beta(s)$ will give rise to a non-trivial extension
\begin{equation}\label{extension}
0 \to \mathcal{O}_C \to \mathcal{E} \to L \to 0.
\end{equation}
Note that $F_*\mathcal{O}_C \otimes L^{-1} \cong F_*(\mathcal{O}_C(-pD))$, and that $\gamma \circ \beta (s) = 0$ by \ref{les1}, so the following sequence will split:
\begin{equation}\label{split}
0 \to \mathcal{O}_C \to F^*\mathcal{E} \to F^*L \to 0.
\end{equation}

\subsection{\texorpdfstring{$\mathbb{P}^1$}{P1}-bundle on \texorpdfstring{$C$}{C}}

Write $P = \mathbb{P}(\mathcal{E})$ for the ruled surface over $C$ with structure map $\pi: P \to C$ and section $\sigma: C \to P$ associated to \ref{extension} with image curve $E = \mathrm{Im}(\sigma) \subset P$. Since the sequence \ref{split} splits, the splitting map $F^*L \to F^*\mathcal{E} \subset \mathrm{Sym}^p(\mathcal{E})$ induces a non-zero section $0 \neq t \in H^0(P, \pi^*L^{-p} \otimes \mathcal{O}_P(p))$ hence a curve $C^{''} \subset P$.

\begin{prop}\cite{R78}\label{standard}
\leavevmode
\begin{itemize}
\item[1.] The two curves $E$ and $C^{''}$ are both smooth, and they are disjoint from each other.
\item[2.] $\omega_{C^{''}/C} = (\pi^*L)|_{C''}$.
\item[3.] $\pi|_{C^{''}}: C^{''} \to C$ is purely inseparable of degree $p$.
\item[4.] $\mathcal{O}_P(E) = \mathcal{O}_P(1)$, and $\mathcal{O}_P(C^{''}) = \mathcal{O}_P(p) \otimes \pi^*L^{-p}$.
\item[5.] $(E^2) = \deg L = \dfrac{2g - 2}{p}$.
\end{itemize}
\end{prop}

\subsection{Cyclic covers over \texorpdfstring{$P$}{P}}
We construct cyclic covers of $P$ that are ramified along the divisor $E + C^{''}$. Suppose $l$ is a positive integer such that $l \mid (p + 1)$ and $l \mid \deg L = \dfrac{2g - 2}{p}$. Write $m = \dfrac{p + 1}{l}$ and let $D \in \abs{L}$ be an effective divisor. Take a line bundle $N \in \mathrm{Pic}(C)$ such that $N^l = L$. Define
\begin{equation}\label{m}
M \coloneqq \mathcal{O}_P(- m) \otimes \pi^*N^{p}.
\end{equation}
Then we have $M^{-l} = \mathcal{O}_P(E + C^{''})$. Define the degree $l$ cyclic cover by $X = \textbf{Spec}(\bigoplus_{i = 0}^{l - 1}M^i) \xrightarrow{\psi} P$. Write $\phi: X \xrightarrow{\psi} P \xrightarrow{\pi} C$ for the composition, and $\tilde{E} = \psi^{-1}(E), C' = \psi^{-1}(C^{''})$ for the reduced preimages of the ramification curves. In particular, $\psi^*(C^{''}) = lC'$ and $\psi^*(E) = l\tilde{E}$.

\begin{prop}\cite{R78, T10}\label{general}
\leavevmode
\begin{itemize}
\item[1.] $X$ is a smooth surface.
\item[2.] $(\tilde{E}^2) = \dfrac{2g - 2}{pl}$.
\item[3.] $\omega_X = \mathcal{O}_X((p - m - 1)l\tilde{E}) \otimes \phi^*\omega_C(- \dfrac{pl - p - l}{l}D)$.
\item[4.] $\omega_X$ is ample if $p \geq 5$ or $(p, l) = (3, 4)$.
\item[5.] $\phi: X \to C$ is a singular fibration with every fiber $F$ having a cuspidal singularity at $F \cap C'$ locally of the form $(X^l = Y^p)$.
\item[6.] The geometric genus of any closed fiber $F$ of $\phi: X \to C$ is $\dfrac{(l - 1)(p - 1)}{2}$.
\end{itemize}
\end{prop}

The following lemma that was crucially used in the argument of Takayama's paper \cite[Lemma 4]{T10} is unfortunately wrong.

\begin{lemma}
For $k \geq 1$, we have
\[
\psi_*\mathcal{O}_X(-k\tilde{E}) = \mathcal{O}_P(-kE) \oplus \left(\bigoplus_{i = 1}^{l - 1}M^i \right).
\]
\end{lemma}

\begin{rmk}
For example, take $k = 2$. If the two sheaves on both sides are isomorphic then they should at least have the same Euler characteristic. Consider the standard short exact sequences:
\begin{align}\label{nonreduced}
0 &  \to \mathcal{O}_{\tilde{E}}(- \tilde{E}) \to \mathcal{O}_{2\tilde{E}} \to \mathcal{O}_{\tilde{E}} \to 0, \\
0 & \to \mathcal{O}_{E}(-E) \to \mathcal{O}_{2E} \to \mathcal{O}_{E} \to 0.
\end{align}
Here we verify that the push forward of the first sequence by $\psi$ is not isomorphic to the second. Since $X$ is smooth $2\tilde{E}$ (and hence $\tilde{E}$) is a Cartier divisor on $X$. By the adjunction formula on $X$:
\begin{equation*}
\chi(\mathcal{O}_{\tilde{E}}(- \tilde{E})) = \deg (-\tilde{E} \mid_{\tilde{E}}) + 1 - p_a(\tilde{E}) = -(g - 1)(1 + \dfrac{2}{pl}).
\end{equation*}
On the other hand,
\[
\chi(\mathcal{O}_{E}(-E)) = -(g - 1)(1 + \dfrac{2}{p}).
\]
Hence $\chi(\mathcal{O}_{2\tilde{E}}) \neq \chi(\mathcal{O}_{2E})$, and this implies $\psi_*\mathcal{O}_{2\tilde{E}}\neq \mathcal{O}_{2E}$. In general $\chi(\mathcal{O}_{k\tilde{E}}) = k(1 - g) - \dfrac{k (k - 1)(g - 1)}{pl}$ and $\chi(\mathcal{O}_{kE}) = k(1 - g) - \dfrac{k(k - 1)(g - 1)}{p}$ for any positive integer $k$.
\end{rmk}

\subsection{Contradiction}Indeed, Theorem 21 in \cite{T10} granting the previous lemma will contradict Proposition \ref{regular}. Following the notations of \cite{T10} suppose $H^1(X, Z_{a, b}^{-1}) \neq 0$. We are going to verify that there exists a collection of numbers $(p, g, l, a, b)$ which results a surface $X$ that satisfies the assumptions of Proposition \ref{regular} yet is not regular. 

\begin{claim}\label{ample}
Take $p = 5, a = 6, b = 3, l = 6, g = 16$, then $Z^4_{6, 3} \otimes \omega_X^{-1} = \mathcal{O}_X(6\tilde{E}) \otimes \phi^*N$ is ample on $X$.
\begin{proof}
In fact, Claim \ref{ample} can be easily verified by the Nakai-Moishezon criterion of ampleness since it intersects both a fiber and a section positively. 
\end{proof}
\end{claim}

According to Proposition \ref{regular}, $X$ will be regular. However, by construction $X$ is a finite cover (of degree $l = 6$ prime to the characteristic $p = 5$) of a ruled surface over a high genus curve ($g = 16$), a regular 1-form should pull-back non-trivially to $X$.

\section{New examples of Kodaira non-vanishing on surfaces}
The correction of \cite[Lemma 4]{T10} is included in the following two lemmas and a proposition.

\begin{lemma}
For any positive integer $k \leq l$
\[
\psi_*\mathcal{O}_{k\tilde{E}} = \bigoplus_{i = 0}^{k}\mathcal{O}_{E}(M^i |_{E}).
\]
In particular,
\[
\psi_*\mathcal{O}_X(-l\tilde{E}) = \bigoplus_{i = 0}^{l - 1}M^i(-E).
\]
\begin{proof}
We proceed the proof with a local computation since the desired identities are of local nature. Suppose $U = \mathrm{Spec}(R) \subset P$ is an open affine subset, and $f \in R$ is the regular function such that $I_{E + C''}|_U = \langle f \rangle$, i.e., $f$ defines the reduced ramification divisor in $U$. Let $t$ be an independent variable such that the inverse image of $U$ in $X$ is
\[
V = \mathrm{Spec}(R[t]/\langle t^l - f \rangle).
\]
We assume that $f$ factors as $f = \bar{s}_1\bar{s}_2$ where $\bar{s}_1$ and $\bar{s}_2$ define $E$ and $C''$ respectively, and that $t$ factors correspondingly as $t = s_1 s_2$ such that $(s_1)_0 = \tilde{E}$ and $(s_2)_0 = C'$.

First the isomorphism $\psi_*\mathcal{O}_{\tilde{E}} = \mathcal{O}_E$ is verified locally on $U$ by the isomorphism $R/\langle \bar{s}_1\rangle \cong R[t]/\langle s_1, t^l - f \rangle$ as $R/\langle \bar{s}_1 \rangle$-modules. We can assume that $l \geq 2$, otherwise the proposition is trivial. On $U$ the coordinate ring of $2\tilde{E}|_U$ is $R[t]/\langle s_1^2, t^l - f \rangle$ and that of $2E|_U$ is $R/\langle \bar{s}_1^2 \rangle$. The ring $R[t]/\langle s_1^2, t^l - f \rangle$ is naturally an $R/\langle \bar{s}_1\rangle$-module. Since $s_1|t$, we have
\[
R[t]/\langle s_1^2, t^l - f \rangle \cong R/\langle \bar{s}_1 \rangle \oplus \langle t \rangle R/\langle \bar{s}_1 \rangle.
\]
Here $f$ is a local generator of the ideal of $E + C''$, so $t$ is a local generator of the invertible sheaf $M$.

If $U_1 = \mathrm{Spec}(R_1) \subset P$ is another open affine subscheme with $V_1 = \mathrm{Spec}(R[t_1]/\langle t_1^l - f_1 \rangle)$, $f_1 = \bar{u}_1\bar{u}_2$, and $t_1 = u_1u_2$ as before, then there exist some unit $\lambda$ such that $t = \lambda t_1$ and $f = \lambda^l f_1$ over the intersection $U \cap U_1$. The transition function for the rank 2 free $R/\langle \bar{s}_1\rangle$-module is the matrix $\mathrm{diag}(1, \lambda)$, which is the same as that of $\mathcal{O}_E \oplus M|_E$.

Hence we prove the statement for the case $k = 2$. For $3 \leq  k \leq l$, it is proved in the same way. The last statement follows from the case $k = l$ and 5-lemma.
\end{proof}
\end{lemma}

\begin{lemma}
For $l + 1 \leq k \leq 2l$,
\[
\psi_*\mathcal{O}_{k\tilde{E}} = \left(\bigoplus_{i = 0}^{k - l - 1}M^i|_{2E} \right)\oplus \left(\bigoplus_{i = k - l}^{l - 1}M^i|_{E} \right).
\]
In general, for any positive integer $k$, write $q = \lfloor \dfrac{k}{l} \rfloor$ and $r = k - ql \geq 0$. Then
\[
\psi_*\mathcal{O}_{k\tilde{E}} = \left(\bigoplus_{i = 0}^{r - 1}M^i|_{(q + 1)E}\right) \oplus \left(\bigoplus_{i = r}^{l - 1}M^i|_{qE}\right).
\]
\begin{proof}
First if $r = 0$, then $k\tilde{E} = \psi^*(qE)$ as a divisor on $X$. Clearly this implies that
\[
\psi_*\mathcal{O}_X(-k \tilde{E}) = \psi_*\mathcal{O}_X(-\psi^*(qE)) = \mathcal{O}_P(-qE) \otimes \left(\bigoplus_{i = 0}^{l - 1}M^i\right).
\]
If $k$ is not divisible by $l$, then we prove the general statement by a local computation as in the previous lemma. With the same notations, note that $R/\langle \bar{s}_1^{a} \rangle$ injects into $R[t]/\langle t^l - f, s_1^k \rangle$ if and only if $a \geq q + 1$. Also note that $R[t]/\langle t^l - f, s_1^k \rangle$ lies in between $R[t]/\langle t^l - f, s_1^{ql} \rangle$ and $R[t]/\langle t^l - f, s_1^{(q + 1)l} \rangle$, which defines the chain of schemes $ql\tilde{E}|_U \subset k\tilde{E}|_U \subset (q + 1)l\tilde{E}|_U$. Since $\psi^*(\bar{s}_1) = s_1^l$, $t^a\bar{s}_1^q$ is divisible by $s_1^k$ if and only if $a + ql \geq k = ql + r$, i.e., $a\geq r$.
\end{proof}
\end{lemma}

\begin{prop}\label{correctlemma}
For any positive integer $k$ write $q = \lfloor \dfrac{k}{l} \rfloor$ and $r = k - ql \geq 0$. Then
\[
\psi_*\mathcal{O}_X(-k\tilde{E}) = \left(\bigoplus_{i = 0}^{r - 1}M^i(- (q + 1)E)\right) \oplus \left(\bigoplus_{i = r}^{l - 1}M^i(- qE) \right).
\]
\begin{proof}
This proposition simply follows from the previous lemmas.
\end{proof}
\end{prop}

\begin{rmk}
Note the contrast between Proposition \ref{correctlemma} and \cite[Lemma 4]{T10}, instead of concentrating at one degree the multiplicities of $\tilde{E}$ will spread as evenly as possible among the summands.
\end{rmk}

\cite[Theorem 18]{T10} is a direct consequence of the erroneous \cite[Lemma 4]{T10}. Correspondingly we correct its statement in the following lemma, which can be proved by a modification of the original proof of \cite[Theorem 18]{T10}.

\begin{lemma}\label{directimage}
Let $Z = \mathcal{O}_X(\tilde{E}) \otimes \phi^*N$ be an ample line bundle on the surface $X$. We write $d = (p + 1)/l$. For any positive integer $n$ we have the form $n = ql + r$ with some non-negative integers $q$ and $r \leq l - 1$. Then
\leavevmode
\begin{itemize}
\item[(1)]
\begin{equation*}
\phi_*Z^{-n} = \left(\bigoplus_{i = 0}^{r - 1} \pi_*\mathcal{O}_P(-(id + q + 1)) \otimes N^{ip - n}\right) \oplus \left( \bigoplus_{k = r}^{l - 1} \pi_*\mathcal{O}_P(- (kd + q)) \otimes N^{kp - n}\right)
\end{equation*}
\item[(2)] For $n \leq l$,
\begin{equation*}
\mathbf{R}^1\pi_*(\psi_*Z^{-n}) = \left(\bigoplus_{i = 1}^{n - 1}\mathrm{Sym}^{id - 1}(\mathcal{E})^{\vee} \otimes N^{ip - n - l}\right) \oplus \left( \bigoplus_{k = n}^{l - 1} \mathrm{Sym}^{kd - 2}(\mathcal{E})^{\vee} \otimes N^{kp - n - l}\right)
\end{equation*}
For $n > l$,
\begin{equation*}
\mathbf{R}^1\pi_*(\psi_*Z^{-n}) = \left(\bigoplus_{i = 0}^{r - 1}\mathrm{Sym}^{id + q - 1}(\mathcal{E})^{\vee} \otimes N^{ip - n - l}\right) \oplus \left( \bigoplus_{k = r}^{l - 1} \mathrm{Sym}^{kd + q - 2}(\mathcal{E})^{\vee} \otimes N^{kp - n - l}\right)
\end{equation*}
\item[(3)]
\begin{equation*}
H^0(C, \mathbf{R}^1\pi_*(\psi_*Z^{-n})) = \left(\bigoplus_{i = 1}^{r - 1} H^0(C, \mathrm{Sym}^{id + q - 1}(\mathcal{E})^{\vee} \otimes N^{ip - n - l})\right) \oplus \left( \bigoplus_{k = r}^{l - 1} H^0(C, \mathrm{Sym}^{kd + q - 2}(\mathcal{E})^{\vee} \otimes N^{kp - n - l})\right)
\end{equation*}
\end{itemize}
\end{lemma}

By the Leray spectral sequence $E_2^{p, q} = H^p(C, \textbf{R}^q\phi_*Z^{-n}) \Rightarrow H^{p + q}(X, Z^{-n})$, there is the following exact sequence
\begin{equation}\label{leray}
0 \to H^1(C, \phi_* Z^{-n}) \to H^1(X, Z^{-n}) \to H^0(C, \textbf{R}^1\phi_*Z^{-n}) \to H^2(C, \phi_*Z^{-n}) = 0
\end{equation}
Using Lemma \ref{directimage} (1), the leftmost term is computed as
\begin{align*}
H^1(C, \phi_*Z^{-n}) = & H^1(C, \pi_*\mathcal{O}_P(- (q + 1)) \otimes N^{-n}) \oplus H^1(C, \pi_*\mathcal{O}_P(-(d + q + 1)) \otimes N^{p - n}) \oplus \dots \\
& \oplus H^1(C, \pi_* \mathcal{O}_P(- ((r - 1)d + q + 1)) \otimes N^{(r - 1)p - n}) \oplus \dots \\
& \oplus H^1(C, \pi_*\mathcal{O}_P(- (rd + q)) \otimes N^{rp - n}) \oplus \dots \\
& \oplus H^1(C, \pi_*\mathcal{O}_P(-((l - 1)d + q )) \otimes N^{(l - 1)p - n}).
\end{align*}
We note this $H^1$ is zero since $\mathrm{Sym}^a(\mathcal{E}) = 0$ for $a < 0$. Now the sequence \ref{leray} reduces to
\[
H^1(X, Z^{-n}) = H^0(C, \textbf{R}^1\phi_*Z^{-n}).
\]
Next note that
\[
\mathbf{R}^1\phi_*Z^{-n} = \mathbf{R}^1\pi_*(\psi_* Z^{-n}).
\]

By Lemma \ref{directimage} (2) and (3) to compute $H^1(X, Z^{-n})$ it suffices to compute the global sections of $\mathbf{R}^1\pi_*(\psi_*Z^{-n})$ on the base curve $C$. Therefore, to claim the non-vanishing of $H^1(X, Z^{-n})$ for some $n$ it suffices to find a non-zero section of $\mathbf{R}^1\pi_*(\psi_*Z^{-n})$. Consequently, we have non-vanishing results for tensor powers of the line bundle $Z$ of a different range than that in \cite[Theorem 19]{T10}:

\begin{theorem}\label{nonvan1}
Under the current construction of $\phi: X \xrightarrow{\psi} P \xrightarrow{\pi} C$, let $Z = \mathcal{O}_X(\tilde{E}) \otimes \phi^*N$. Then $Z$ is ample. Moreover, for any positive integer $n \leq \lfloor \dfrac{l}{2} \rfloor$, we have
\[
H^1(X, Z^{-n}) \neq 0.
\]
\begin{proof}
By the numerical nature of ampleness we note that the line bundle $Z = \mathcal{O}_X(\tilde{E}) \otimes \phi^*N$ is ample. By construction, $ld = p + 1$ and $n = lq + r$. In the case that $q = 0$ and $n = r \leq \lfloor\dfrac{l}{2} \rfloor$, we have
\[
l((l - n)d - 2) = (l - n)p - n - l,
 \]
so there is a surjection
\[
\mathrm{Sym}^{(l - n)d - 2}(\mathcal{E}) \to N^{(l - n)p - n - l} \to 0.
\]
Dualizing it and twisting it with $N^{(l - n)p - n - l}$ we get a non-zero section:
\[
H^0(C, \mathrm{Sym}^{(l - n)d - 2}(\mathcal{E})^{\vee} \otimes N^{(l - n)p - n - l}) \neq 0.
\]
\end{proof}
\end{theorem}

Next we investigate the bundles $Z_{a, b} \coloneqq \mathcal{O}_X(a\tilde{E}) \otimes \phi^*N^b$. The $Z_{a, b}$ is still ample as long as $a$ and $b$ are positive integers. We aim to find some $a, b$ so that $H^1(X, Z_{a, b}^{-1}) \neq 0$ by standard Leray spectral sequence method again (cf. \cite[Theorem 21]{T10}).
\begin{theorem}\label{nonvan2}
For any pairs of positive integers $a$ and $b$ such that $a \leq l - 1$ and $b \leq l - a$, the ample line bundle $Z_{a, b}$ on the surface $X$ provides a counterexample to the Kodaira vanishing theorem, i.e.,
\[
H^1(X, Z^{-1}_{a, b}) \neq 0.
\]
\begin{proof}
We do the division with remainder for $a$ with respect to $l$ as before: $a = ql + r$, for some non-negative integers $q$ and $r \leq l - 1$. Also write $d = (p + 1)/l$. Then
\begin{align*}
\mathbf{R}^1\phi_*Z^{-1}_{a, b} & = \mathbf{R}^1\pi_*(\psi_*Z^{-1}_{a, b}) = \mathbf{R}^1\pi_*(\psi_*\mathcal{O}_X(-a\tilde{E})) \otimes N^{-b}\\
& = \left(\bigoplus_{i = 0}^{r - 1}\mathbf{R}^1\pi_*\mathcal{O}_P(- id - q - 1) \otimes N^{ip - b}\right) \oplus \left(\bigoplus_{i = r}^{l - 1}\mathbf{R}^1\pi_*\mathcal{O}_P(- id - q) \otimes N^{ip - b} \right).
\end{align*}
The computation of global sections is similar as before:
\begin{align*}
H^0(C, \mathbf{R}^1\pi_*(\psi_*Z^{-1}_{a, b})) & = H^0(C, \mathrm{Sym}^{d + q - 1}(\mathcal{E})^{\vee} \otimes N^{p - b - l}) \oplus \dots \oplus H^0(C, \mathrm{Sym}^{(r - 1)d + q - 1}(\mathcal{E})^{\vee} \otimes N^{(r - 1)p - b - l}) \\
& \oplus H^0(C, \mathrm{Sym}^{rd + q - 2}(\mathcal{E})^{\vee} \otimes N^{rp - b - l}) \oplus \dots \oplus H^0(C, \mathrm{Sym}^{(l - 1)d + q - 2}(\mathcal{E})^{\vee}\otimes N^{(l - 1)p - b - l}),
\end{align*}
and we have that $H^0(C, \mathbf{R}^1\pi_*(\psi_*Z^{-1}_{a, b})) \neq 0$ for $a \leq l - 1$ and $b \leq l - a$.
\end{proof}
\end{theorem}

\begin{rmk}
A direct calculation shows that
\[
Z^{p - 1}_{a, b} \otimes\omega_X^{-1} = \mathcal{O}_X((ap - a - pl + p + l + 1)\tilde{E}) \otimes \phi^*\mathcal{O}_C((bp - b + pl - p - l)N - K_C).
\]
We can easily verify that the multiplicity $ap - a - pl + p + l + 1$ of $\tilde{E}$ and the degree of $\mathcal{O}_C((bp - b + pl - p - l)N - K_C)$ cannot be both positive under the range of the preceding theorem. This is compatible with Proposition \ref{regular}.
\end{rmk}

The rest of this section is to show that the direct image of the relative dualizing sheaf $\phi_*(\omega_{X/C}^k)$ of the fibration $\phi: X \to C$ is not nef for any positive integer $k$.

\begin{cor}\label{nonnef}
Let $\phi: X \to C$ be the singular fibration as before. Then $\phi_*(\omega_{X/C}^{k})$ is not nef for any positive integer $k$.
\begin{proof}
The computation is the same for all $k > 0$. Below we exemplify the case that $k = 1$ for simplicity of notations. By part 3 of Proposition \ref{general},
\[
\phi_*\omega_{X/C} = \pi_*(\psi_*\mathcal{O}_X((p - m - 1)l\tilde{E})) \otimes \mathcal{O}_C(- \dfrac{pl - p - l}{l}D).
\]
By projection formula we note
\[
\psi_*\mathcal{O}_X(l\tilde{E})) = \psi_*(\psi^*\mathcal{O}_P(E)) = \bigoplus_{i = 0}^{l - 1}(M^i \otimes \mathcal{O}_P(E)).
\]
By the definition of $M$ (\ref{m}) we obtain
\begin{equation}\label{decomp}
\phi_*\omega_{X/C} = \bigoplus_{i = 1}^{l} \left(\mathrm{Sym}^{p - im - 1}(\mathcal{E}) \otimes N^{ip + l - pl} \right).
\end{equation}
Taking $i = 1$, we focus on one particular direct summand $V \coloneqq \mathrm{Sym}^{p - m - 1}(\mathcal{E}) \otimes N^{p + l - pl}$ of the decomposition \ref{decomp}. We have isomorphism $\mathbb{P}(V) \cong \mathbb{P}(\mathrm{Sym}^{p - m - 1}(\mathcal{E}))$ with $\mathcal{O}_{\mathbb{P}(V)}(1)$ corresponding to $\mathcal{O}_{\mathbb{P}(\mathrm{Sym}^{p - m - 1}(\mathcal{E}))}(1) \otimes \pi_{p - m - 1}^*(N^{p + l - pl})$, where $\pi_{p - m - 1}: \mathbb{P}(\mathrm{Sym}^{p - m - 1}(\mathcal{E})) \to C$ is the usual projective bundle map. In particular, $\pi_1: \mathbb{P}(\mathcal{E}) \to C$ coincides with the previous notation $\pi: P \to C$, which is compatible over $C$ with the $(p - m - 1)$-th Veronese embedding $\iota_{p - m - 1}: \mathbb{P}(\mathcal{E}) \to \mathbb{P}(\mathrm{Sym}^{p - m - 1}(\mathcal{E}))$. Restricting $\mathcal{O}_{\mathbb{P}(\mathrm{Sym}^{p - m - 1}(\mathcal{E}))}(1) \otimes \pi_{p - m - 1}^*(N^{p + l - pl})$ to the image of $\iota_{p - m - 1}$ we obtain the line bundle
\[
W \coloneqq \mathcal{O}_{\mathbb{P}(\mathcal{E})}(p - m - 1) \otimes \pi^*(N^{p + l - pl}).
\]
Intersecting $W$ with the section $E$ we have
\[
(W \cdot E) = (p - m - 1)(E^2) + (p + l - pl) \cdot \dfrac{\deg L}{l},
\]
where by part 5 of Proposition \ref{standard} $(E^2) = \deg L$. Hence
\[
(W \cdot E) = \deg L \cdot \left[(p - m - 1) + \dfrac{p + l - pl}{l}\right] = -\dfrac{\deg L}{l} < 0.
\]
This shows that $V$ as a quotient of $\phi_*\omega_{X/C}$ is not nef. Hence neither is $\phi_*\omega_{X/C}$.\end{proof}
\end{cor}

\begin{rmk}
The data $(C, D \in \abs{L}, f \in K(C))$ is a 1-dimensional example of a \textit{TR-triple} defined for arbitrary dimension in \cite[Sec. 2]{Mu13}. In fact, Mukai constructs an $(n + 1)$-dimensional TR-triple based on an $n$-dimensional TR-triple in an inductive manner under certain assumption on the exact 1-form $df$, generalizing the Tango-Raynaud structure for curves and surfaces (see \cite[Prop. 1.7, Prop. 2.3]{Mu13}). The statement and the proof of Corollary \ref{nonnef} remain valid if the triple $(C, D \in \abs{L}, f \in K(C))$ is replaced by any higher dimensional TR-triple.
\end{rmk}

As a corollary to Corollary \ref{nonnef}, we show that our surface also provides a counterexample of Koll\'ar's vanishing theorem. Note that
\begin{equation}\label{n}
L^{p - m - 1} \otimes N^{p + l - pl} = N^{l(p - m - 1) + p + l - pl} = N^{-1}.
\end{equation}

\begin{cor}\label{kollar}
Notations as before. Then
\[
H^1(C, \phi_*\omega_X \otimes N) \neq 0.
\]
\begin{proof}
By \ref{decomp} and \ref{n}, $N^{-1}$ is a quotient of $V$, which in turn is a quotient of $\phi_*\omega_{X/C}$. Hence $\omega_C \cong (\omega_C \otimes N^{-1}) \otimes N$ is a quotient of $V \otimes \omega_C \otimes N$, which is in turn a quotient of $\phi_*\omega_X \otimes N$. The nonvanishing of $H^1(C, \omega_C)$ implies the desired result.
\end{proof}
\end{cor}

\section{Comments and remarks}

In this section we review a few ideas behind the various geometric constraints of the Kodaira non-vanishing and deduce some direct consequences which are at most implicit in \cite{E88, ShB91a, Mu13}. For that matter, we present the results all in the form of corollaries.

\subsection{Torsors of infinitesimal group schemes and foliations}

Ekedahl studied the pluri-canonical systems on surfaces generalizing Bombieri's results on their base-point freeness and very ampleness, and considered vanishing of cohomologies of line bundles on surfaces. Suppose $X$ is a minimal surface of general type and $L$ is an ample line bundle on $X$ such that $H^1(X, L^{-1}) \neq 0$. The following sequence is exact in the flat topology:
\begin{equation}\label{torsor}
0 \to \alpha_{L^{-1}} \to L^{-1} \to L^{-p} \to 0,
\end{equation}
where $\alpha_{L^{-1}}$ is an infinitesimal group scheme of order $p$ over $X$. Ekedahl analyzed the case where $H^1(X, L^{-p}) = 0$ and $H^1(X, L^{-1}) \neq 0$. The induced long exact sequence in cohomology of \ref{torsor} implies that $H^1(X, \alpha_{L^{-1}}) \neq 0$. The group $H^1(X, \alpha_{L^{-1}}) \neq 0$ classifies torsors of $\alpha_p$ over $X$, which are purely inseparable covers over $X$ of degree $p$. Such a purely inseparable degree $p$ map $\pi: Y \to X$ from a normal surface $Y$ corresponds to a $1$-foliation $\mathcal{F}$ on $X$, which induces a quotient map $\phi: X \to Y' \coloneqq X/\mathcal{F}$. Then $Y$ is birational to $Y'^{(1)}$, the pre-image of $Y'$ under the Frobenius map. For details see \cite{E88}.

The existence of an ample line bundle as a $1$-foliation on the surface $X$, under some mild conditions on $X$ and the foliation, could result $X$ to be (purely inseparably) uniruled. As a corollary to Ekedahl's results, we have the following effective vanishing for surfaces not of general type:

\begin{cor}\label{effective}
Suppose that $X$ is a smooth projective surface not of general type, and that $L$ is an ample line bundle on $X$. Then
\[
H^1(X, L^{-2}) = 0.
\]
\begin{proof}
By \cite[Ch. II, Thm 1.6]{E88}, the only remaining case to check is that $X$ is quasi-elliptic, $\kappa(X) = 1$, and $p = 2$ or 3. Let $\phi: X \to C$ be a quasi-elliptic fibration, and $F$ be a general fiber of $\phi$, which is a rational curve with a single ordinary cusp. Then $(F^2) = 0$, $(F \cdot K_X) = 0$ and $p_a(F) = 1$.

Since $L$ is ample on $X$, $H^1(X, L^{-n}) = 0$ for all sufficiently large $n$ by Serre vanishing and Serre duality. If the first cohomologies of all negative tensor powers of $L$ vanish, then there is nothing to prove. Now we assume that there exists an integer $e \geq 0$ such that $H^1(X, L^{-n}) \neq 0$ and $H^1(X, L^{-pn}) = 0$. By \cite[Ch. II, Thm 1.3 (v)]{E88}, we have the following inequality:
\begin{equation}\label{inequality}
-2 \leq (K_X \cdot F) + (1 - p)(L^{n} \cdot F) + (F^2) = (1 - p)(L^{n} \cdot F)
\end{equation}

If $p = 3$, then there is no solution for \ref{inequality} if $n \geq 2$.

If $p = 2$, then \ref{inequality} holds true only for $(L \cdot F) = 1$ and $n = 2$. Suppose that $H^1(X, L^{-2}) \neq 0$. Let $Y$ be the normalization of the relative Frobenius fiber product $X \times_C C^{(1)}$ over the Frobenius map $F_C: C^{(1)} \to C$ on $C$. We write $F_{X/C}: Y \to X$ for the generically finite degree 2 purely inseparable morphism. Since the singularities on the fibers of $f$ are ordinary cusps, $g: Y \to C^{(1)}$ is a geometrically ruled surface. Write $\mathcal{F}$ for the rank 2 vector bundle on $C^{(1)}$ associated to $g$ and $A \in \mathrm{Pic}(C^{(1)})$ for the divisor corresponding to $\mathrm{det}(\mathcal{F})$. In particular,
\begin{equation}\label{canonicalruled}
K_Y \sim -2C_0 + g^*(K_{C^{(1)}} + A),
\end{equation}
where $C_0$ is a fixed section of $g$ with $\mathcal{O}_{\mathbb{P}(\mathcal{F})}(1) \cong \mathcal{O}_Y(C_0)$.

On the other hand, $Y$ is birational to an $\alpha_2$-torsor $\mu: Z \to X$ induced by a non-zero class in $H^1(X, L^{-2})$. Suppose $\alpha: Z' \to Z$ is a minimal resolution of singularities of $Z$, and denote the composition $Z' \to Z \to C^{(1)}$ by $h$, and the purely inseparable morphism $Z' \to X$ by $\tau$.

By \cite[Ch. I, Prop. 2.1]{E88},
\begin{equation}\label{canonicaltorsor}
2K_{Z'} \sim \tau^*(K_X - 2L) - B
\end{equation}
for some curve $B$ on $Z'$ exceptional for $\alpha$. By the canonical bundle formula for quasi-elliptic surfaces (see \cite[Thm. 2]{BM77}) we have
\begin{equation}\label{canonicalqell}
K_X \sim f^*(K_C - M) + \sum_{\lambda}P_{\lambda},
\end{equation}
where $M = \sheafhom_{\mathcal{O}_C}(f_*\omega_X, \mathcal{O}_C) \otimes \omega_C$ (by duality, see \cite[Formula 2.2.3]{DR73}) is the torsion-free part of $\mathbf{R}^1f_*\mathcal{O}_X$, $\sum_{\lambda}P_{\lambda}$ is supported in the multiple fibers of $f$ (for a priori, the multiplicity of each connected component $P_{\lambda}$ in $\sum_{\lambda}P_{\lambda}$ is less than its multiplicity as a multiple fiber, which is equal to 2 by \cite[Cor. 5.1.2]{CD89}, in particular, every multiple fiber $2P_{\lambda}$ with $2P_{\lambda}$ appearing in \ref{canonicalqell} is tame). Plugging \ref{canonicalqell} into \ref{canonicaltorsor} and using the commutativity of $f \circ F_{X/C} = F_C \circ g: Y \to C$, we can compare the canonical bundles of $Y$ and $Z'$. Let $W$ be another birationally ruled surface birationally dominating both $Y$ and $Z$ with maps $h_1: W \to Z'$ and $g_1: W \to Y$.

Then by \ref{canonicalruled} and \ref{canonicaltorsor} we have
\begin{equation}\label{comparison}
h_1^*\tau^*(2L) - 4g_1^*(C_0) \sim -2(g \circ g_1)^*(A) - 2(h \circ h_1)^*(M) - h_1^*(B) + (\tau \circ h_1)^*(\sum_{\lambda}P_{\lambda}) \mod \mathrm{Exc}(h_1, g_1),
\end{equation}
where $\mathrm{Exc}(h_1, g_1)$ is the collection of exceptional curves of $g_1$ and $h_1$ in $W$ and \ref{comparison} should be understood as a linear equivalence possibly up to a divisor that is supported in $\mathrm{Exc}(h_1, g_1)$. Since $C$ is a smooth curve of genus at least $1$, every irreducible exceptional curve in $\mathrm{Exc}(h_1, g_1)$ is contained in some fiber of one of the compositions $h \circ h_1$ and $g \circ g_1$, hence the right hand side of \ref{comparison} consists of vertical divisors. Looking at the left hand side of \ref{comparison}, the horizontal component of $h_1^*\tau^*(2L)$ must be the same as $4g_1^*(C_0)$. Moreover, every component in the vertical part of $h_1^*\tau^*(2L)$ is effective since $L$ is ample, which implies the same for the right hand side. Next notice that the proof of Corollary \ref{nonnef} works in this case ($p = 2, l = 3, m = 1$), showing that $M$ is nef (more precisely, it has degree $1$ on $C^{(1)}$). Also the divisor $A$ is ample since $(C_0^2) > 0$, being the reduced preimage of a section in $X$. Thus \ref{comparison} is reduced to
\begin{equation}
\mathrm{Ver}(h_1^*\tau^*(2L)) + 2(g \circ g_1)^*(A) + 2(h \circ h_1)^*(M) \sim (\tau \circ h_1)^*(\sum_{\lambda}P_{\lambda}) \mod \mathrm{Exc}(h_1, g_1, \alpha),
\end{equation}
where $\mathrm{Ver}(h_1^*\tau^*(2L))$ consists of entire fibers in $h_1^*\tau^*(2L)$. In particular, the number of components on the right is even. By \ref{canonicalqell}, the intersection number of $K_X$ with a section will be odd. However, by \ref{canonicaltorsor} (or in fact, before passing to the resolution, $2K_Z \sim \tau^*(K_X - 2L)$, by abuse of notation still denote the covering map by $\tau$), the degree of $K_X$ on a section is even, hence a contradiction.
\end{proof}
\end{cor}

\begin{rmk}
This improves part of \cite[Cor. 5.9]{DiCF15} for the case that the surface is not of general type. After finishing this argument, Adrian Langer informed the author about a stronger statement in his coming paper with a simpler proof.
\end{rmk}

\subsection{Stability of rank 2 vector bundles}

A non-zero class $\eta \in H^1(X, L^{-1})$ gives rise to an extension of locally free sheaves
\[
0 \to \mathcal{O}_X \to \mathcal{E} \to L \to 0.
\]
One computes that $\delta(\mathcal{E}) \coloneqq c_1^2(\mathcal{E}) - 4c_2(\mathcal{E}) = (L^2) > 0$. By \cite[Theorem 1]{ShB91a}, there exists a positive integer $n$ and an integral surface $Y \subset \mathbb{P}(F^{n*}\mathcal{E})$ such that the restriction $\phi: Y \to X$ is purely inseparable of degree $p^n$. By Riemann-Roch, Shepherd-Barron shows that (see \cite[Prop. 14]{ShB91a})
\begin{equation}\label{insepcan}
\chi(\mathcal{O}_Y) = p^n \chi(\mathcal{O}_X) + \dfrac{1}{12}p^n(p^n - 1)[(2p^n - 1)(L^2) - 3(L \cdot K_X)].
\end{equation}
\begin{prop}\cite[Prop. 19]{ShB91a}
\[
K_Y \equiv \phi^*(K_X + (1 - p^n)L).
\]
\end{prop}

\begin{cor}
Suppose $X$ is a minimal surface of general type and $(p^k - 1)L - K_X$ is nef and big for some positive integer $k$. Then $H^1(X, L^{-1}) = 0$.
\begin{proof}
Suppose to the contrast. By taking $n$ sufficiently large, the assumption guarantees that one can always assume $Y$ to be a ruled surface. In particular, $H^0(Y, K_Y) = 0$, and consequently $\chi(\mathcal{O}_Y) \leq 1$. On the other hand, we want to show that $\chi(\mathcal{O}_Y) \geq 2$ using \ref{insepcan}.

\textit{Case 1}. $p \geq 3$. By \cite[Cor. 1.5]{G14}, we have $\chi(\mathcal{O}_X) > 0$. Then it is elementary to check that \ref{insepcan} yields that $\chi(\mathcal{O}_Y) \geq 2$.

\textit{Case 2}. $p = 2$. By \cite[Ch. II Cor. 1.8]{E88}, we have $\chi(\mathcal{O}_X) \geq 2 - K_X^2$. Plugging-in \ref{insepcan} we have $\chi(\mathcal{O}_Y) \geq 2$ as well for $n$ sufficiently large.
\end{proof}
\end{cor}

\subsection{Singular fibration structure}
Suppose $D \in \abs{L}$ is an effective divisor. Then there exists $f \in K(X)$ such that $(df) > pD$, where $(df)$ is the divisor associated to the exact differential form. The rational function $f$ induces, after blowing up the loci of indeterminacy, a morphism $X' \to \mathbb{P}^1$. The Stein factorization of this morphism gives rise to a fibration $\phi: X' \to C$, where $X'$ is birational to $X$. Mukai showed that the generic fiber is singular, as the sheaf $\Omega^1_{X'}/\phi^*(\Omega^1_C)$ has torsion.

\begin{rmk}
Suppose the divisor $D$ is locally defined by $s = 0$ for $s \in H^0(X, L)$, then the function $s^{p + 1}$ will suffice to locally define the desired exact 1-form, as $d(s^{p + 1}) = (p + 1)s^p(ds)$.

The map induced by such $s$ is only rational in priori, but thanks to the following statement due to Mukai (see \cite[Prop. 3.2]{Mu13}): if $f: X' \to X$ is a blowup of smooth surfaces, then vanishing of $H^1$ for any big and nef line bundles on $X$ is equivalent to the vanishing of $H^1$ for any big and nef line bundles on $X'$.
\end{rmk}

In fact, if $H^1(X, L) = 0$ for any big and nef line bundle $L$ on $X$, then $H^1(X', f^*L) = 0$ and $f^*L$ is big and nef on $X'$. So the implication from $X'$ to $X$ is clear. The proof of the other direction uses the sheaf $B^1$ and that
\[
H^0(X, B^1 \otimes L^{-1}) = \{df \mid f \in K(X), (df) \geq pD\}
\]
for $D \in \abs{L}$.

\begin{rmk}
Generic singular fibrations first appeared in Bombieri-Mumford's classification of quasi-elliptic surfaces, where the singularities are all ordinary cusps. In general, the singularities appearing in the general fiber of such fibration could only be \textit{unibranched}. Suppose in general $(p^k - 1)L - K_X$ is pseudo-effective for some positive integer $k$. The singularities in the fiber of $\phi: X' \to C$ can be resolved by finitely many iterations of the relative Frobenius morphism over the base curve as in the proof of \cite[Lemma 9]{ShB91b}. In particular, $X$ is uniruled. This \textit{genus-change} phenomenon was first observed by Tate (see \cite{T52}) and reproved using the language of schemes by Schr\"{o}er (see \cite{S09}).
\end{rmk}

\begin{cor}
Suppose $\abs{L} \neq \emptyset$, $(p^k - 1)L - K_X$ is pseudo-effective for some integer $k > 0$ and $0 \neq \eta \in \mathrm{ker}[H^1(X, L^{-1}) \to H^1(X, L^{-p})]$ is a non-zero class giving rise to a fibration $\phi: X' \to C$. Suppose $\phi_* \mathcal{O}_{X'} \cong \mathcal{O}_C$. Then $C$ is a Tango curve.
\begin{proof}
Suppose $D \in \abs{L}$ is an effective divisor on $X$. By assumption, there exist a rational function $f \in K(X)$ such that $(df) \geq pD$ and a normal birationally ruled surface $\tau: Y \to C$ with a purely inseparable morphism $h: X' \to Y$ of degree $p^k$, such that the $k$-th iteration of the Frobenius map of $X$ factors through $h$ and $\phi = \tau \circ h$. There is an integral curve $E \subset X'$ corresponding to the singularities on the fibers of $\phi$ such that $\phi |_E: E \to C$ is purely inseparable of degree $p^k$. Write $E'$ for the reduced image of $E$ in $Y$. Then $f^*(E') = E$, $f_*(E) = p^k E'$, and $E'$ is a section of $\sigma$. Note that locally the sheaf $\phi^*\Omega^1_{C}$ is generated by $df$, in particular $f \in K(C)$. If $D$ contains a component that is supported on a section $\sigma: D_1 \to X'$ of $\phi$, then we have $(df)|_{D_1} \geq L|_{D_1}$. Hence $D_1$ and also $C$ is a Tango curve. If $D$ does not contain such a component, then $E \subset D$. In the same vein we see that $E$ is a Tango curve. So $C$ is a Tango curve.
\end{proof}
\end{cor}

\subsection{Global vector fields and differential 1-forms}\label{globalvec}

The current construction is compatible with that of Takeda (see \cite{T92}). Using exactly the same proof as those of \cite[Thm. 2.1, Cor. 3.4]{T92}, the surface $X$ is an example of a smooth surface of general type with \textit{non-trivial vector fields} and \textit{non-closed global differential 1-forms} (we refer the reader to the original paper of Takeda for details).

\begin{rmk}
To the author's knowledge, almost all known examples of Kodaira non-vanishing live in isolation, so it would be very interesting to have such examples in families. It might not be easy to find such in the moduli of surfaces, but the seemingly more accessible question is to find examples of Kodaira non-vanishing on a fixed surface with the pathological polarization varying in a positive dimensional family in the Picard group.
\end{rmk}

% ----------------------------------------------------------------
\bibliographystyle{amsalpha}

\end{document}